\documentclass[reqno,12pt]{amsart}

\usepackage{graphicx, amsmath, amssymb, amscd, amsthm, euscript, psfrag, amsfonts,bm}

\newtheorem{thm}[equation]{Theorem}

\newtheorem{Con}[equation]{Conjecture}
\newtheorem{rmk}[equation]{Remark}
\newtheorem{cor}[equation]{Corollary}
\newtheorem{prop}[equation]{Proposition}

\numberwithin{equation}{section}

\newcommand{\be}{begin{equation}}

\newcommand{\bH}{\mathbb H}

\newcommand{\q}{\mathbb{Q}}

\newcommand{\e}{\epsilon}
\newcommand{\z}{\mathbb{Z}}
\renewcommand{\q}{\mathbb{Q}}

\newcommand{\br}{\mathbb{R}}

\newcommand{{\grinv}}{{\Cal G}^{-r}}

\newcommand{\ba}{\backslash}

\newcommand{\G}{\Gamma}

\newcommand{\Cal}{\mathcal}

\newcommand{\la}{\langle}
\newcommand{\ra}{\rangle}

\newcommand{\SL}{\operatorname{SL}}
\newcommand{\GL}{\operatorname{GL}}
\newcommand{\bp}{\begin{pmatrix}}
\newcommand{\ep}{\end{pmatrix}}
\renewcommand{\be}{\begin{equation*}}
\newcommand{\ee}{\end{equation*}}

\renewcommand{\bp}{{\rm bp}}

\newcommand{\SO}{\operatorname{SO}}

\newcommand{\op}{\operatorname}

\newcommand{\gq}{\Gamma(q)}\newcommand{\Gq}{\Gamma(q)}
\begin{document}

\title[Multiplicity]{Eigenvalues of congruence covers of geometrically finite hyperbolic manifolds}

\author{Hee Oh}
\address{Mathematics department, Brown university, Providence, RI 02912
and Korea Institute for Advanced Study, Seoul, Korea}
\email{heeoh@math.brown.edu}

\begin{abstract} Let $G=\SO(n,1)^\circ$ for $n\ge 2$ and $\G$ a geometrically finite Zariski dense subgroup of $G$
which is contained in an arithmetic subgroup of $G$. 
Denoting by $\G(q)$ the principal congruence subgroup of $\G$ of level $q$, and fixing a positive number $\lambda_0$ strictly smaller than $(n-1)^2/4$,
we show that, as $q\to \infty$ along primes,
the number of Laplacian eigenvalues of the congruence cover $\G(q)\ba \bH^n$ smaller than $\lambda_0$
is at most of order $[\G:\Gq]^c$ for some $c=c(\lambda_0)>0$. 
\end{abstract}

\maketitle
\section{Introduction}
Let $G$ denote identity component of the special orthogonal group $\SO(n,1)$ for $n\ge 2$.
As well-known, $G$ is the group of orientation preserving isometries of the real hyperbolic space $\bH^n$.
Denote by $\Delta$ the negative of the Laplace-Beltrami operator of $\bH^n$. 
On any complete hyperbolic manifold
$M$ of dimension $n$,
$\Delta$ acts on the space of smooth functions on $M$ with compact support and admits a unique extension
to an unbounded self-adjoint positive operator on $L^2(M)$. We denote by $\sigma(M)$ its spectrum of $M$. For instance,
$\sigma(\bH^n)$ is known to be $[\tfrac{(n-1)^2}4, \infty)$.


There exists a torsion-free discrete subgroup $\G$ of $G$ such that $M=\G\ba \bH^n$. The limit set
$\Lambda(\G)$ of $\G$ is the smallest non-empty closed $\G$-invariant subset of the geometric boundary of $\bH^n$.
The convex core $\mathcal C(M)$ of $M$ is 
the quotient by $\G$ of the smallest convex subset of $\bH^n$ containing all geodesics
connecting points in $\Lambda(\G)$.
We say that $M$ (or $\G$) is {\it geometrically finite} 
if the unit neighborhood of $\mathcal C(M)$ has finite volume.
Geometrically finite manifolds are natural generalizations of manifolds with finite volume. 

For a geometrically finite hyperbolic manifold $M$ of infinite volume, Lax and Phillips \cite{LP} showed that
$\sigma(M)$ is the disjoint union of the discrete spectrum consisting of finitely many eigenvalues
contained in $[0,\tfrac{(n-1)^2}4)$ with finite multiplicities
and the essential spectrum which is a closed sub-interval of $[\tfrac{(n-1)^2}4, \infty)$.

In this note, we are interested in the number of discrete eigenvalues of a sequence of congruence covers of a fixed hyperbolic manifold, which
itself is a covering of an arithmetic hyperbolic manifold of finite volume.

Let $G$ be defined over $\q$ via a $\q$-embedding $G\to \GL_N$ for some positive integer $N$. We set $G(\z)=G\cap \GL_N(\z)$. Fixing a subgroup $\G$ of $G(\z)$, and a positive integer $q$,
we consider the $q$-th principal congruence subgroup $\G(q)$ of $\G$ given by
$$\G(q):=\{\gamma\in \G: \gamma \equiv e \mod q\}.$$

\bigskip
We denote by $\delta=\delta(\G)$ the critical exponent of $\G$, i.e., the abscissa of convergence of the Poincar\'e series
$\mathcal P(t)=\sum_{\gamma\in \G}e^{-t d(o, \gamma(o))}$ for $o\in \mathbb H^n$. 
In the rest of the introduction, we assume that $\G$ is torsion-free, geometrically finite and Zariski dense in $G$.
For a fixed $0<\lambda_0 <\tfrac{(n-1)^2}4$,
denote by $$\mathcal N (\lambda_0, \Gq)$$ 
the number of
eigenvalues of $\Gq\ba \bH^n$ contained in the interval $[0,\lambda_0]$, counted with multiplicities.
Here is our main theorem:
\begin{thm}\label{main} 
There exists $\eta>0$ such that for any $0<\lambda_0<\frac{(n-1)^2}4$,
\begin{equation}\label{t1} \mathcal N(\lambda_0, \Gq) \ll [\G : \gq]^{\frac{\delta -s_0}{\eta}} \quad\text{for any $q$ prime} ,\end{equation}
where $\tfrac{n-1}{2}<s_0\le n-1$ such that $\lambda_0=s_0(n-1-s_0)$ and the implied constant depends only on $\lambda_0$.
\end{thm}

In \cite{Ha}, Hamenst\"adt showed that the number of eigenvalues of $\G\ba \bH^n$ 
smaller than $\lambda_0$
is bounded from above by $c^{\text{Vol}(\mathcal C (\G\ba \bH^n))}$
for some $c>0$ depending only on the dimension $n$.
In case when $\G$ is a lattice, a stronger upper bound of $c\cdot \text{Vol}(\G\ba \bH^n)$
was previously known by Buser-Colbois-Dodziuk \cite{BCD}.

Since $\Gq$ is a normal subgroup of $\G$ of finite index, the limit set of $\Gq$ is equal to the limit set
of $\G$, and hence
$$\frac{\text{Vol}(\mathcal C(\Gq\ba \bH^n))}{\text{Vol}(\mathcal C (\G\ba \bH^n))}= [\G:\Gq] .$$
Therefore Theorem \ref{main} implies the following stronger upper bound for congruence coverings of
an arithmetic manifold $\G\ba \bH^n$ of infinite volume: as $q\to \infty$ along primes,
\begin{equation}
 \mathcal N(\lambda_0, \Gq) \ll  \op{Vol}(\mathcal C (\Gq\ba \bH^n))^{\frac{\delta -s_0}{\eta}}.\end{equation}

\begin{rmk}\rm
\begin{enumerate}
 \item 
We can relax the restriction on $q$ so that
$q$ is square-free with no divisors from a fixed finite set of primes.

\item
As $\Gq$ is geometrically finite, the discrete spectrum of $\Gq\ba \bH^n$ is non-empty only when $\delta>\tfrac{n-1}{2}$,
in which case, $\delta(n-1-\delta)$ is the smallest discrete eigenvalue and has multiplicity one \cite{Su}.
\item In our proof of Theorem \ref{main}, $\eta$ can be taken to be any number smaller than 
the uniform spectral gap $$\liminf_{q:\text{primes}} (\delta-s_{1q})$$
where $s_{1q}<\delta$ is such that $s_{1q}(n-1-s_{1q})$ is the second smallest eigenvalue of $\Gq\ba \bH^n$.
The existence of a {\it positive} uniform spectral gap follows from
the works of Bourgain-Gamburd \cite{BG} 
and of Bourgain-Gamburd-Sarnak \cite{BGS2} for $n=2$.
 Their result has been generalized by Salehi-Golsefidy-Varju \cite{SV} for a general connected semisimple algebraic group.
These methods do not provide an explicit estimate on $\eta$. However
for $\delta$ large,  Gamburd \cite{Ga} ($n=2$) and Magee \cite{Ma} $(n\ge 3$) gave a lower bound on $\eta$ when $\G$ is contained in certain arithmetic subgroups of $G$.
For instance, according to \cite{Ga},
if $\G$ is contained in $\SL_2(\z)$ and $\delta>\tfrac 56$, then $\eta>\delta-\tfrac{5}{6}$.
\end{enumerate}
\end{rmk}

Note that any finitely generated subgroup $\G$ of $\SL_2(\br)$ with $\delta>0$ is geometrically finite,
and Zariski dense in $\SL_2(\br)$.
If $\G$ is contained in $\SL_2(\z)$, the strong approximation property of Zariski dense subgroups (\cite{MV} and \cite{No})
implies that $\G(q)\ba \G$ is isomorphic to $\SL_2(\mathbb F_q)$ for all but finitely many primes $p$.
Hence $[\G:\Gq]=\# \SL_2(\mathbb F_q)\ll q^3$.

Therefore we deduce from Theorem \ref{main} and the above remark 1.4.(3):
\begin{cor} Let $\G<\SL_2(\z)$ be a torsion-free finitely generated subgroup with $\delta>\tfrac 56$.
Then for any $1/2< s_0\le \delta$,
$$  \mathcal N(s_0(1-s_0), \Gq) \ll q^{\tfrac{18(\delta-s_0)}{6\delta -5}}\quad \text{ for any prime $q$}.$$
\end{cor}

For the following discussion, we find it convenient to put $\rho=\tfrac{n-1}{2}$ and to use
the parametrization $\lambda_s=s(n-1-s)$ so that
$$[0, \tfrac{(n-1)^2}4)=\{ \lambda_s : s\in (\rho,n-1]\}.$$
For each $s\in (\rho,n-1]$, we denote by $$m(s, \G(q))$$ the multiplicity of
$\lambda_s$ occurring as a discrete eigenvalue of $\Delta$ in $\Gq\ba \bH^n$.
Note that $$\mathcal N(s_0(n-1-s_0),\Gq)=\sum_{s_0\le s\le n-1} m(s,\Gq).$$
In the influential paper of Sarnak and Xue \cite{SX},
a conjecture for the upper bound of $m(s,\G(q))$ is formulated for arithmetic groups of a connected
 semisimple algebraic group $G$ defined over $\q$.
In our setting of $G=\SO(n,1)^\circ$, their conjecture can be stated as follows:
\begin{Con}[Sarnak-Xue] \label{up} If $[G(\z):\G]<\infty$, then for any fixed $s\in (\rho,n-1]$,
as $q\to \infty$,
$$ m(s,\gq)\ll_{\e} [\G : \gq]^{\frac{(n-1)-s}{\rho}+\e}  \quad \text{for any $\e>0$} .$$\end{Con}

Fixing $o\in \bH^n$ and denoting by $d$ the hyperbolic distance in $\bH^n$,
consider the lattice point counting function:
 $$N(T, q):=\#\{\gamma\in \Gq: d(\gamma(o), o)\le T \}.$$ 

\begin{Con}[Sarnak-Xue] \label{up2} Let $[G(\z):\G]<\infty$. For any $T\gg 1$ and $q\gg 1$,
$$ N(T,q) \ll_\e \frac{e^{(n-1+\e)T}}{[\G:\Gq ]} +e^{\rho T}\quad \text{for any $\e>0$} $$
where the implied constant is independent of $T$ and $q$.
\end{Con}

Sarnak and Xue proved that if $\G$ is a uniform lattice in $G$, then
Conjecture \ref{up2} implies
Conjecture \ref{up}; in fact, their methods prove a stronger statement that
the same upper bound works for $\mathcal N(s(n-1-s),\Gq)$ as well.
This implication has been extended to all lattices in $G$ by Huntley and Katznelson \cite{HK}. 
Sarnak and Xue \cite{SX} proved Conjecture \ref{up2} when $\G$ is a cocompact arithmetic subgroup of $\SO(n,1)$ for $n=2,3$ and
hence settled Conjecture \ref{up} in this case.


In view of the conjectures \eqref{up} and \eqref{up2} and the above discussion,
we pose the following two conjectures: let $\G<G(\z)$ be geometrically finite and Zariski dense in $G$.
\begin{Con} \label{con1} For any $s\in (\rho,\delta]$,
\begin{equation*} \mathcal N(s(n-1-s),\gq)\ll_{\e} [\G : \gq]^{\frac{\delta -s}{\delta-\rho}+\e}  \quad \text{for any $\e>0$}.
            \end{equation*} \end{Con}

\begin{Con}\label{con2}  For $T\gg 1$,
\begin{equation*} N(T,q) \ll_\e \frac{e^{(\delta +\e)T}}{[\G:\Gq ]} +e^{\rho T}\quad \text{for any $\e>0$.} 
            \end{equation*} \end{Con}

\begin{prop}\label{f} 
Conjecture \ref{con2} implies Conjecture \ref{con1}.
\end{prop}

Note that Conjecture \ref{con1} implies the following limit formula, which is also suggested by the Plancherel formula of $L^2(G)$ given by Harish-Chandra:
\begin{Con} \label{con3} For any $s\in (\rho,\delta]$,
\begin{equation*} \lim_{q\to \infty}\frac{ \mathcal N(s(n-1-s),\gq)}{ [\G : \gq]} =0.
            \end{equation*} \end{Con}
Conjecture \ref{con3} is known to be true if $\G$ is a co-compact lattice by DeGeorge and Wallach \cite{DW}
or if $\G=\SL_2(\z)$ by Sarnak \cite{Sar}.
It does not seem to be known for a general (arithmetic) lattice, although
Savin proved it for those $s$ whose corresponding eigenfunctions are cusp-forms \cite{Sa}.
We also refer to \cite{FLW} where an analogous problem was answered positively for $\G=\SL_n(\z)$.

\bigskip
The proof of Conjecture \ref{up} by Sarnak-Xue \cite{SX} for co-compact arithmetic lattices of $\SO(2,1)$
and $\SO(3,1)$ uses number theoretic arguments which give a very sharp uniform upper bound (Conjecture \ref{up2})
for the number of lattice points $\Gq$ in a ball, using explicit realizations of arithmetic groups
as units of certain division algebras over number fields. The presence of $\rho$ in the denominator
of the exponent in Conjecture \ref{up} (Theorem for the cases in discussion) is due to the sharpness of this counting technique.

This approach won't be possible for general Zariski dense subgroups, as such an explicit arithmetic realization is not available
for this rather wild class of groups (which makes one wonder that Conjecture \ref{con2} is perhaps too bold).

Instead, we use a recent work of Mohammadi and the author \cite{MO}
where uniform counting results for orbits of $\Gq$'s were obtained with an error term.

Another important ingredient is the so-called {\it Collar Lemma} on uniform estimates on the size
of eigenfunctions away from flares and cusps of $\Gq\ba \bH^n$;
this was first obtained by Gamburd
in \cite{Ga} for $\SO(2,1)$ and generalized by Magee \cite{Ma} for all $\SO(n,1)$.

\section{}
In the rest of this paper, let $\G$ be geometrically finite and Zariski dense, with $\delta>(n-1)/2$.
We continue notations from the introduction.

\subsection{Lattice point counts}
As before, set $$N(T,q):=\# \G(q)\cap B_T$$
where $B_T:=\{g\in G: d(g(o),o)\le T\}$.

We first recall the following lattice point counting theorem in \cite{MO}:
\begin{thm}\label{mo}  There exists $\eta>0$ such that for any prime $q$,
\be N(T,q)= c\cdot \frac{e^{\delta T}}{[\G:\Gq]} +O({e^{(\delta-\eta) T}})
 \quad \text{for all $T\gg 1$} \ee
where both $c>0$ and the implied constant are independent of $q$.
\end{thm}

We remark that the above type of lattice point counting theorem is stated
in \cite{MO} under a uniform spectral gap hypothesis. However
since $B_T$ is a bi-$K$-invariant subset, one needs only a uniform {\it spherical}
spectral gap in proving Theorem \ref{mo}. And the uniform spherical spectral gap
property for $L^2(\Gq\ba \bH^n)$ for $q$ primes follows from \cite{SV} by the transfer
principle from combinatorial spectral gap to an archimedean one, obtained by Bourgain, Gamburd and Sarnak 
\cite{BGS2} (also see \cite{Kim} for $n\ge 3$).

\subsection{} 
We identify $\bH^n=G/K$ for $K=\SO(n)$ and denote by $o\in \bH^n$ whose stabilizer is $K$.
Fix a Haar measure $dg$ on $G$. This induces an invariant measure on $\Gq\ba G$ for which we use
the same notation $dg$ by abuse of notation. The right translation action of $G$ on $L^2(\Gq\ba G, dg)$
preserves the measure $dg$, and hence yields a unitary representation of $G$.
Let $\mathcal C$ denote the Casimir operator of $G$. The action of $\mathcal C$ on $K$-invariant smooth functions on $G$ is 
given by $-\Delta$. 
For each $s\in (\rho, \delta]$, we denote
by $\pi_s$ the spherical complementary series representation of $G$ on which $\mathcal C$ acts by the scalar $-\lambda_s$. Then
the multiplicity $m(s,\Gq)$ is equal to the multiplicity of $\pi_s$ occurring as a sub-representation of $L^2(\Gq\ba G)$.
This follows from the well-known correspondence between positive definite spherical functions of $G$ and
the spherical unitary dual of $G$.

Define the bi-$K$-invariant function of $G$: $$\psi_s(g):=\la \pi_s(g)v_s, v_s\ra$$ where $v_s$ is the unique $K$-fixed unit vector
in $\pi_s$, up to a scalar.
Then $\psi_s (g) $ is a positive function such that for all $g\in G$,
\begin{equation}\label{es}\psi_s(g) \asymp e^{(s-2\rho)d(o, g(o))}\end{equation}
(we could not find a reference for this
supposedly well-known fact except for \cite{MO}).

\subsection{Proof of Theorem \ref{main}} 
We need to show that
there exists $\eta>0$ such that for any $\rho<s_0\le n-1$,
we have, as $q\to \infty$ along primes,
\be  \sum_{s\in [s_0, n-1]} m(s,\gq)\ll  [\G : \gq]^{\frac{\delta -s_0}{\eta}}  ,\ee
where the implied constant is independent of $q$.

We follow a general strategy of \cite{SX} (also see \cite{Ga}).

Consider the following bi-$K$-invariant functions of $G$:
 $$f_o(g):=\chi_{B_T}(g)\psi_{s_0}(g)\quad \text{and}\quad F_o(g):=f_o*\check{f_o} (g)$$ where $\chi_{B_T}$ is the characteristic function of $B_T$
and $\check f_o(g):=\overline{f_o(g^{-1})}$.

By \cite{SX}, we have
\begin{equation}\label{sx} F_o(g)\ll \begin{cases} e^{2(s_0-\rho)T} e^{-\rho d(g(o),o)} &\text{ if $d(o, g(o))\le 2T$}\\
0&\text{ if $d(o, g(o))> 2T$.}                   \end{cases}\end{equation}

The spherical transform $\hat f$ of a bi-$K$-invariant function $f$ is given by
$$\hat f(\lambda_s)= \int_G f(g) \psi_{s}(g) dg.$$
Hence $$\hat f_o(\lambda_s)=
\int_G \chi_{B_T}(g) \psi_{s_0}(g)\psi_{s}(g) dg $$
and the associated spherical transform
$\hat F_o(\lambda_s)$ is given by $|\hat f_o(\lambda_s)|^2$.

By \eqref{es},
it follows that for all $s\ge s_0$,
\begin{equation}\label{fs}\hat F_o(\lambda_s)\gg  e^{(2s_0+2s -4\rho)T}\ge  e^{(4s_0 -4\rho)T} \end{equation}
where the implied constant depends only on $s_0$.

Define the automorphic kernel $K_q$ on $\Gq\ba G \times \Gq \ba G$ as follows:
$$K_q(g_1, g_2):=\sum_{\gamma\in \Gq} F_o(g_1^{-1}\gamma g_2).$$
Note $K_q(g_1k_1, g_2k_2)=K_q(g_1, g_2)$ for any $k_1, k_2\in K$.

Let $\{\lambda_{j,q} \}$ be the multi-set of discrete eigenvalues of $\Gq\ba\bH^n$ which is finite by
Lax-Phillips and let
$\{\phi_{j,q}\}$ be corresponding real-valued eigenfunctions with $L^2$-norm one in $L^2(\Gq\ba \bH^n)$.
We may understand $\phi_{j,q}$ as a function on $\Gq\ba G$ which is right $K$-invariant.
Let $s_{j,q}\in (\rho, \delta]$ be such that $\lambda_{j,q}=s_{j,q}(n-1-s_{j,q})$.

A key technical ingredient we need is the following result of Gamburd for $n=2$ \cite{Ga}
and of Magee for $n$ general \cite{Ma}:
\begin{thm}\label{ma}
Fix a closed interval $I\subset (\rho, \delta]$.
There exists a compact subset $\Omega\subset \G\ba G$ and $C>0$ such that
for any integer $q\ge 1$ and any $s_{j,q}\in I$,
$$\int_{\Omega_q} |\phi_{j,q}(g)|^2 dg\ge C$$
where $\Omega_q:=\pi_q^{-1}(\Omega)$ for the canonical projection $\pi_q:\Gq\ba G\to \G\ba G$.
\end{thm}

By applying the pretrace formula to $K_q$, we deduce 
$$K_q(g, g)= \sum_{ \lambda_{j,q}} \hat F_o(\lambda_{j,q}) |\phi_{j,q}(g)|^2 +\mathcal E$$
where the term $\mathcal E$ is the contribution from the continuous spectrum.
The positivity of $\hat F_o$ yields that $\mathcal E\ge 0$, and hence
$$K_q(g, g)\ge  \sum_{ \lambda_{j,q}} \hat F_o(\lambda_{j,q}) |\phi_{j,q}(g)|^2 .$$
By integrating over the compact subset $\Omega_q$ of
Theorem \ref{ma}, we obtain
$$\int_{\Omega_q} K_q(g,g) d g\ge 
C\cdot \sum_{ \lambda_{j,q}, s_{j,q}\in I} \hat F_{o}(\lambda_{j,q}) .$$

Therefore setting $I:=[s_0, \delta]$,
we deduce from \eqref{fs} the following:
$$\int_{\Omega_q} K_q(g,g) dg \gg   e^{4(s_{0} -\rho)T}\left( \sum_{ s\in I}  
m(s,\Gq) \right) .$$

On the other hand, using \eqref{sx},
\begin{align*}
\int_{\Omega_q} K_q(g,g) dg&=\sum_{\gamma\in \Gq}\int_{\Omega_q} F_o(g^{-1}\gamma g)dg
\\ &=\sum_{\gamma\in \Gq}\sum_{\gamma_0\in\Gq\ba \G}\int_{\Omega} F_o(g^{-1}\gamma_0^{-1}\gamma \gamma_0 
g)dg \\ &=
[\G:\Gq]\cdot
\sum_{\gamma\in \Gq} \int_{\Omega} F_o(g^{-1}\gamma 
g)dg \\
&\le [\G:\Gq]\cdot e^{2(s_0 -\rho)T} \int_{x\in \Omega, d(x,\gamma x)\le 2T} e^{-\rho d(x,\gamma x)}dx .
 \end{align*}
Using the compactness of $\Omega$, and Theorem \ref{mo},
we deduce, for some $R>0$ (depending only on $\Omega$),
\begin{align}\label{fc} &\int_{x\in \Omega, d(x,\gamma x)\le 2T} e^{-\rho d(x,\gamma x)}dx
\ll\int_{0}^{2T+R}e^{-\rho t} N( T,q) dt
\notag \\ & \ll \tfrac{1}{[\G:\Gq]}e^{(2\delta-2\rho)T} +e^{(2\delta-2\rho-2\eta)T}
\end{align}

Putting these together, we have
$$\left( \sum_{ s\ge s_0}  
m(s,\Gq) \right) \ll e^{(-2s_0+2\delta)T} +{[\G:\Gq]}\cdot e^{(-2s_0+2\delta-2\eta)T}.$$
Hence by setting $T$ so that
$e^{2T}={[\G:\Gq]}^{1/\eta}$, we deduce
$$\left( \sum_{ s\ge s_0}  
m(s,\Gq) \right)\ll {[\G:\Gq]}^{\frac{\delta-s_0}{\eta}},$$
as desired, proving Theorem \ref{main}.

Note that replacing the use of Theorem \ref{mo} by Conjecture \ref{con2} for an upper bound of $N(T,q)$ in \eqref{fc} 
yields $$\mathcal N(s_0(n-1-s_0),\Gq)\ll_\e {[\G:\Gq]}^{\frac{\delta-s_0}{\delta-\rho} +\e},$$
which is Conjecture \ref{con1}.
This verifies Proposition \ref{f}.

\end{document}